# MINIMAX ESTIMATION OF LINEAR FUNCTIONALS OVER NONCONVEX PARAMETER SPACES[1]

By T. Tony Cai and Mark G. Low

*University of Pennsylvania*

The minimax theory for estimating linear functionals is extended to the case of a finite union of convex parameter spaces. Upper and lower bounds for the minimax risk can still be described in terms of a modulus of continuity. However in contrast to the theory for convex parameter spaces rate optimal procedures are often required to be nonlinear. A construction of such nonlinear procedures is given. The results developed in this paper have important applications to the theory of adaptation.

**1. Introduction.** Let $Y$ be an observation from either the white noise model,

$$dY(t) = f(t)\,dt + n^{-1/2}\,dW(t) \tag{1}$$

where $W(t)$ is a standard Brownian motion, or the Gaussian sequence model

$$Y(i) = f(i) + n^{-1/2}\varepsilon_i \tag{2}$$

where $\varepsilon_i$ are i.i.d. standard normal random variables.

The minimax theory for estimating a linear functional $T$ has been studied in great generality when it is assumed that the function $f$ belongs to a parameter space which is convex. See, for example, Ibragimov and Has'minskii (1984), Donoho and Liu (1991a, b) and Donoho (1994). In particular, the properties of the minimax linear estimators can often be described precisely. In this case for any linear functional $T$ write $R_A^*(n; \mathcal{F})$ for the minimum (over all linear procedures) maximum mean squared error. Donoho and Liu (1991a) introduced a modulus of continuity

$$\omega(\varepsilon, \mathcal{F}) = \sup\{|T(g) - T(f)| : \|g - f\|_2 \leq \varepsilon, f \in \mathcal{F}, g \in \mathcal{F}\}$$

Received February 2002; revised December 2002.
[1]Supported in part by NSF Grant DMS-03-06576.

where the norm in this equation is the $L_2$ norm in function space for the white noise with drift model and the $l_2$ norm in sequence space for the sequence model. Donoho and Liu (1991a, b) and Donoho (1994) have shown that in either of these two cases,

$$(3) \qquad R_A^*(n; \mathcal{F}) = \sup_{\varepsilon > 0} \omega^2(\varepsilon, \mathcal{F}) \frac{1/(4n)}{1/n + \varepsilon^2/4}$$

and that

$$\frac{1}{8}\omega^2\left(\frac{2}{\sqrt{n}}, \mathcal{F}\right) \leq R_A^*(n; \mathcal{F}) \leq \frac{1}{4}\omega^2\left(\frac{2}{\sqrt{n}}, \mathcal{F}\right).$$

An earlier version of this result can also be found in Ibragimov and Has'minskii (1984). Without the restriction to affine procedures write $R_N^*(n; \mathcal{F})$ for the minimax mean squared error for estimating the linear functional $T$. Donoho and Liu (1991b) have shown that

$$\frac{R_A^*(n; \mathcal{F})}{R_N^*(n; \mathcal{F})} \leq 1.25.$$

Therefore the maximum risk of the optimal linear procedure is within a small constant factor of the minimax risk when the parameter space is convex. Of equal importance, Donoho and Liu (1991b) showed that the modulus can be used to give a recipe for constructing an affine procedure which has the maximum mean squared error attaining the bound given in (3).

Recent work on estimating linear functionals has focused on adaptive estimation. The goal is to find a single procedure which is near minimax simultaneously over a number of different parameter spaces. Pioneering work in this area began with Lepski (1990). This work focused on particularly important examples such as Lipschitz classes. In Efromovich and Low (1994) a general theory was developed for the case of nested convex parameter spaces.

A general extension of this adaptive estimation theory to spaces which are not nested must also include a minimax analysis for sets which are not convex. The reason for this is that we need to first know the minimax risk over the union of the original convex spaces and this space need not be convex unless the sets are nested. This paper focuses on such an extension of the minimax theory for estimating linear functionals over nonconvex parameter spaces. For applications to adaptive estimation see Cai and Low (2002). Although as just mentioned our primary motivation for this problem is the theory of adaptation the minimax theory itself is in fact quite interesting. In particular in this setting optimal linear procedures can sometimes have risks far from the optimal rate. In fact even if the parameter space is only a union of two convex sets it is possible that the maximum risk of the best linear estimator does not even converge even though the maximum risk of the optimal nonlinear procedure converges quickly. Such examples are given in Section 5.

Although optimal linear procedures need no longer be close to optimal we show that the minimax rate of convergence is still determined by the modulus of continuity over the parameter space when the parameter space is a finite union of convex sets. On the other hand, in Section 4, it is shown that the minimax linear risk is determined by the modulus of continuity over the convex hull of the parameter space. Therefore affine procedures fail when, in terms of the modulus, the convex hull is much larger than the parameter space itself. Such are the cases in the examples in Section 5. In these cases rate optimal estimators need to be nonlinear. A general construction of such nonlinear procedures is given in Section 3.

One of the main tools for the construction of the general procedure is a construction of linear procedures which have a given variance and precisely control the bias over two different convex parameter spaces. Upper bounds are given on the bias over one parameter space and lower bounds over the other. These linear procedures can then be used to test which of the convex sets the function lies in and then usual linear procedures can be used. The details of these arguments can be found in Sections 2 and 3.

The theoretical results are complemented by several illustrative examples given in Section 5 covering a range of cases. In the examples of estimating a linear functional of a nearly black object the parameter space is the union of a growing number of convex parameter spaces. In these cases the usual minimax lower bound is no longer sharp and the minimax rate of convergence is derived explicitly using a mixture prior and a constrained risk inequality.

**2. Ordered modulus and bias variance tradeoffs.** One of the main tools for the construction of the general minimax procedure is the construction of linear procedures which have a given variance and precisely control the bias over two different convex parameter spaces. Upper bounds are given on the bias over one parameter space and lower bounds over the other. The key technical tool which allows for this construction is an ordered modulus of continuity between two function spaces. It is a generalization of the modulus of continuity introduced by Donoho and Liu (1991a) which has already been shown in Low (1995) to allow for the construction of a procedure which minimizes the maximum squared bias given a constraint on the maximum variance.

For a linear functional $T$ define an ordered modulus of continuity between two classes $\omega(\varepsilon, \mathcal{F}, \mathcal{G})$ by

$$\omega(\varepsilon, \mathcal{F}, \mathcal{G}) = \sup\{Tg - Tf : \|g - f\|_2 \leq \varepsilon; f \in \mathcal{F}, g \in \mathcal{G}\}.$$

Note that $\omega(\varepsilon, \mathcal{F}, \mathcal{G})$ does not necessarily equal $\omega(\varepsilon, \mathcal{G}, \mathcal{F})$. It is clear that the modulus $\omega(\varepsilon, \mathcal{F}, \mathcal{G})$ is an increasing function of $\varepsilon$ and $0 \leq \omega(\varepsilon, \mathcal{F}, \mathcal{G}) \leq \infty$ if $\mathcal{F} \cap \mathcal{G} \neq \varnothing$. The between class modulus is also instrumental in the analysis of adaptation over different parameter spaces [see Cai and Low (2002)].

When $\mathcal{G} = \mathcal{F}$, $\omega(\varepsilon, \mathcal{F}, \mathcal{F})$ is the usual modulus of continuity over $\mathcal{F}$ and will be denoted by $\omega(\varepsilon, \mathcal{F})$. The following result on the concavity of the modulus is important in the bias variance tradeoffs and in the construction of the general minimax procedure.

THEOREM 1. *Assume that $\mathcal{F}$, $\mathcal{G}$ are convex and that $\mathcal{F} \cap \mathcal{G} \neq \varnothing$. Let $T$ be a linear functional. Then the function $\omega(\varepsilon, \mathcal{F}, \mathcal{G})$ is a concave function of $\varepsilon$. In particular it follows that, for $D > 1$,*

$$\omega(D\varepsilon, \mathcal{F}, \mathcal{G}) \leq D\omega(\varepsilon, \mathcal{F}, \mathcal{G}).$$

PROOF. Suppose that $g_1 \in \mathcal{G}$, $g_2 \in \mathcal{G}$ and $f_1 \in \mathcal{F}$, $f_2 \in \mathcal{F}$ with

$$\|g_i - f_i\|_2 \leq \varepsilon_i.$$

Then, for $0 \leq \lambda \leq 1$,

$$\|\lambda g_2 + (1-\lambda)g_1 - [\lambda f_2 + (1-\lambda)f_1]\|_2 \leq \lambda \varepsilon_2 + (1-\lambda)\varepsilon_1$$

and

$$T(\lambda g_2 + (1-\lambda)g_1 - [\lambda f_2 + (1-\lambda)f_1])$$
$$= \lambda(T(g_2) - T(f_2)) + (1-\lambda)(T(g_1) - T(f_1)).$$

It then follows that

$$\omega(\lambda \varepsilon_2 + (1-\lambda)\varepsilon_1, \mathcal{F}, \mathcal{G}) \geq \lambda \omega(\varepsilon_2, \mathcal{F}, \mathcal{G}) + (1-\lambda)\omega(\varepsilon_1, \mathcal{F}, \mathcal{G})$$

and so $\omega$ is concave.  □

As mentioned earlier in Low (1995) it was shown that in the white noise model for any linear functional the modulus of continuity can be used to precisely trade off various levels of bias and variance over a given convex parameter space. The modulus of continuity between parameter spaces can be used to perform an analogous trade. It can be used to give a linear procedure which has upper bounds for the bias over one parameter space and lower bounds for the bias over the other parameter space. The detailed results are given in Theorems 2 and 3 below.

We shall write $\langle u, v \rangle$ for the usual $l_2$ inner product for either sequence or function space. Specifically if we observe the white noise with drift model let

$$\langle f, g \rangle = \int fg$$

and if we observe the sequence model let

$$\langle f, g \rangle = \sum f_i g_i.$$

For all $V \geq 0$ let

(4) $$B(V, \mathcal{F}, \mathcal{G}) = 2^{-1} \sup_{\varepsilon > 0} (\omega(\varepsilon, \mathcal{F}, \mathcal{G}) - \sqrt{nV}\varepsilon).$$

It will also be convenient to introduce an inverse of $B(V, \mathcal{F}, \mathcal{G})$ defined for all $B \geq 0$ by

(5) $$V(B, \mathcal{F}, \mathcal{G}) = \sup_{\varepsilon > 0} \frac{1}{n\varepsilon^2} ([\omega(\varepsilon, \mathcal{F}, \mathcal{G}) - 2B]_+)^2.$$

We shall show in Theorems 2 and 3 that there is a linear estimator with variance bounded by $V$, which has maximum bias over $\mathcal{F}$ less than or equal to $B(V, \mathcal{F}, \mathcal{G})$ and minimum bias over $\mathcal{G}$ greater than or equal to $-B(V, \mathcal{F}, \mathcal{G})$. Theorem 2 covers the most usual situations where linear estimators can be easily described in terms of the modulus. Theorem 3 extends the theory to cover the general case.

Our analysis is split into a number of cases. The most usual ones are covered by cases 1(a) and 2(a). It is these cases which are in fact needed in the construction of the general procedure in Section 3. We include the others for completeness. First note that we shall always assume that $\omega(1, \mathcal{F}, \mathcal{G}) > 0$; otherwise the linear functional is constant over $\mathcal{F} \cup \mathcal{G}$ and the estimation problem is thus trivial.

CASE 1. Suppose that $0 < B(V, \mathcal{F}, \mathcal{G}) < \infty$. Then define $\varepsilon(V, \mathcal{F}, \mathcal{G})$ by

(6) $$\varepsilon(V, \mathcal{F}, \mathcal{G}) = \arg\max_{\varepsilon \geq 0} (\omega(\varepsilon, \mathcal{F}, \mathcal{G}) - \sqrt{nV}\varepsilon)$$

where $\varepsilon(V, \mathcal{F}, \mathcal{G})$ is the smallest value of $\varepsilon$ for which the maximum in (4) is attained. It will be convenient to break case 1 into two further cases, namely:

(a) $0 < \varepsilon(V, \mathcal{F}, \mathcal{G}) < \infty$.
(b) $\varepsilon(V, \mathcal{F}, \mathcal{G}) = \infty$.

CASE 2. $B(V, \mathcal{F}, \mathcal{G}) = 0$ and $B(V', \mathcal{F}, \mathcal{G}) > 0$ for all $0 \leq V' < V$. Note that if $B(V', \mathcal{F}, \mathcal{G}) = 0$ for some $V' < V$ then we could reduce the variance of our estimator without increasing the magnitude of the bias. Under this assumption there are only two possibilities.

(a) $\omega(\varepsilon, \mathcal{F}, \mathcal{G}) = \sqrt{nV}\varepsilon$ on some interval $0 \leq \varepsilon \leq \varepsilon_0$ where $\varepsilon_0 > 0$. We can then define $\varepsilon(V, \mathcal{F}, \mathcal{G})$ to be the largest $\varepsilon \leq \frac{1}{\sqrt{n}}$ for which $\omega(\varepsilon, \mathcal{F}, \mathcal{G}) = \sqrt{nV}\varepsilon$.

(b) $\omega(\varepsilon,\mathcal{F},\mathcal{G}) < \sqrt{nV}\varepsilon$ whenever $\varepsilon > 0$. It then follows from the concavity of the modulus that $0 < B(V',\mathcal{F},\mathcal{G}) < \infty$ for some $V' < V$. In this case set $\varepsilon(V,\mathcal{F},\mathcal{G}) = 0$.

The following technical lemma shows that $B(V,\mathcal{F},\mathcal{G})$ is continuous in $V$ whenever it is finite.

LEMMA 1. *Suppose $\mathcal{F}$ and $\mathcal{G}$ are closed and convex with $\mathcal{F} \cap \mathcal{G} \neq \varnothing$. Then $\varepsilon(V,\mathcal{F},\mathcal{G})$ is nonincreasing in $V$. Assume $B(V,\mathcal{F},\mathcal{G}) < \infty$. Then $\varepsilon(V', \mathcal{F},\mathcal{G}) < \infty$ if $V' > V$ and*

$$\lim_{V_m \downarrow V} B(V_m,\mathcal{F},\mathcal{G}) = B(V,\mathcal{F},\mathcal{G}). \tag{7}$$

*If, in addition, $B(V',\mathcal{F},\mathcal{G}) < \infty$ for some $V' < V$, then*

$$\lim_{V_m \to V} B(V_m,\mathcal{F},\mathcal{G}) = B(V,\mathcal{F},\mathcal{G}). \tag{8}$$

PROOF. Note first that the monotonicity of $\varepsilon(V,\mathcal{F},\mathcal{G})$ and the fact that $\varepsilon(V',\mathcal{F},\mathcal{G}) < \infty$ if $V' > V$ follows from the concavity of the modulus $\omega(\varepsilon,\mathcal{F},\mathcal{G})$ as shown in Theorem 1. Now assume that $B(V,\mathcal{F},\mathcal{G}) < \infty$ and let $V_m \downarrow V$. Note that

$$B(V_m,\mathcal{F},\mathcal{G}) \leq B(V,\mathcal{F},\mathcal{G}) \qquad \text{for } V_m \geq V,$$

and that for any $\varepsilon$

$$B(V_m,\mathcal{F},\mathcal{G}) \geq 2^{-1}(\omega(\varepsilon,\mathcal{F},\mathcal{G}) - \sqrt{nV_m}\varepsilon).$$

Taking limits yields

$$\liminf_{V_m \downarrow V} B(V_m,\mathcal{F},\mathcal{G}) \geq 2^{-1}(\omega(\varepsilon,\mathcal{F},\mathcal{G}) - \sqrt{nV}\varepsilon)$$

for all $\varepsilon$ and so taking the supremum over all $\varepsilon$ on the right-hand side shows that the limit exists and is equal to $B(V,\mathcal{F},\mathcal{G})$. This proves (7).

Note that $B(V,\mathcal{F},\mathcal{G})$ is a convex function of $\sqrt{V}$ since it is a supremum of a collection of convex functions of $\sqrt{V}$. Hence $B(V,\mathcal{F},\mathcal{G})$ is continuous in $V$ on any open interval over which it is finite. Hence if $B(V',\mathcal{F},\mathcal{G}) < \infty$ for some $V' < V$ then $B(\cdot,\mathcal{F},\mathcal{G})$ is continuous at $V$ and so (8) follows. □

We now state the bias—variance tradeoff theorem in the most easily understood and most typical case where $0 < \varepsilon(V,\mathcal{F},\mathcal{G}) < \infty$ and the modulus is attained by two functions $f \in \mathcal{F}$ and $g \in \mathcal{G}$.

THEOREM 2. *Suppose $\mathcal{F}$ and $\mathcal{G}$ are convex and closed with $\mathcal{F} \cap \mathcal{G} \neq \varnothing$. Assume that $0 < \varepsilon(V,\mathcal{F},\mathcal{G}) < \infty$. Suppose further that there are $f \in \mathcal{F}, g \in \mathcal{G}$ such that*

$$\|g - f\|_2 = \varepsilon(V,\mathcal{F},\mathcal{G}) \equiv \varepsilon_V \quad \text{and} \quad Tg - Tf = \omega(\varepsilon_V,\mathcal{F},\mathcal{G}). \tag{9}$$

*Write $u \equiv \frac{g-f}{\varepsilon_V}$ for the direction of the affine family joining $g$ and $f$. Let*

$$a = T\left(\frac{f+g}{2}\right) - \sqrt{nV}\left\langle u, \frac{f+g}{2}\right\rangle. \tag{10}$$

*Then the estimator*

$$\hat{T}_V = a + \sqrt{nV}\int u(t)\,dY(t) \tag{11}$$

*for the white noise with drift model and the estimator*

$$\hat{T}_V = a + \sqrt{nV}\sum u(i)Y(i) \tag{12}$$

*for the sequence model have constant variance*

$$E(\hat{T}_V - E\hat{T}_V)^2 = V \tag{13}$$

*and have biases bounded by*

$$\sup_{f\in\mathcal{F}} E_f \hat{T}_V - Tf = B(V,\mathcal{F},\mathcal{G}) \tag{14}$$

*and*

$$\inf_{g\in\mathcal{G}} E_g \hat{T}_V - Tg = -B(V,\mathcal{F},\mathcal{G}). \tag{15}$$

REMARK. If $\mathcal{F}$ and $\mathcal{G}$ are closed, convex and norm bounded with nonempty intersection then the condition that the modulus is attained is guaranteed. The extension to cases where either the modulus is not attained as well as for when $\varepsilon(V,\mathcal{F},\mathcal{G}) = 0$ and $\varepsilon(V,\mathcal{F},\mathcal{G}) = \infty$ will be covered in Theorem 3.

PROOF OF THEOREM 2. The proof of this theorem essentially follows that of Theorem 2 in Low (1995). Note that the proofs of (14) and (15) are entirely similar so we shall only give the details for the proof of (15).

Let $f \in \mathcal{F}$ and $g \in \mathcal{G}$ be extremal functions satisfying (9) which exist since $\mathcal{F}$ and $\mathcal{G}$ are closed. Let $h$ be any other element of $\mathcal{G}$. The affine family joining $g$ and $h$ is given by $(1-\theta)g + \theta h$, $0 \leq \theta \leq 1$. Let

$$J(\theta) = T((1-\theta)g + \theta h) - Tf - \sqrt{nV}\|(1-\theta)g + \theta h - f\|_2.$$

It follows from the definition of $\varepsilon(V,\mathcal{F},\mathcal{G})$ given in (6) that $J(\theta) \leq J(0)$ for all $0 \leq \theta \leq 1$ and since $J(\theta)$ is clearly differentiable it follows that $J'(0) \leq 0$. A simple computation shows that

$$Th - Tg - \sqrt{nV}\langle u, (h-g)\rangle \leq 0. \tag{16}$$

Now

$$E\hat{T}_V - Tg = T\left(\frac{f+g}{2}\right) + \sqrt{nV}\left\langle u, \left(g - \frac{f+g}{2}\right)\right\rangle - Tg \tag{17}$$

and

(18) $$E\hat{T}_V - Th = T\left(\frac{f+g}{2}\right) + \sqrt{nV}\left\langle u, \left(h - \frac{f+g}{2}\right)\right\rangle - Th.$$

It then follows from (16)–(18) that

(19) $$(E\hat{T}_V - Tg) - (E\hat{T}_V - Th) \leq 0.$$

Finally note that a simple calculation yields

(20) $$E\hat{T}_V - Tg = -B(V, \mathcal{F}, \mathcal{G}).$$

Equations (19) and (20) combine to show (15) and the proof is complete. □

Theorem 2 treats the cases 1(a) and 2(a) under the additional assumption that the modulus is attained by $f \in \mathcal{F}$ and $g \in \mathcal{G}$. The functions $f$ and $g$ are used explicitly in the construction of the estimate $\hat{T}_V$. In general, the modulus may not be attained and in these cases the description of a linear estimator which trades variance and bias is more involved. We describe the general case in detail in the following theorem. Some of the details are similar to those given in Section 12 of Donoho (1994).

Define $B(m)$ to be the closed $L_2$ ball with radius $m$ and let $\mathcal{F}_m = \mathcal{F} \cap B(m)$ and $\mathcal{G}_m = \mathcal{G} \cap B(m)$. It follows from Lemma 2 of Donoho (1994) that for $\mathcal{F}_m$ and $\mathcal{G}_m$ the modulus $\omega(\varepsilon, \mathcal{F}_m, \mathcal{G}_m)$ can always be attained by some $f \in \mathcal{F}_m$ and $g \in \mathcal{G}_m$.

Define $V_m$, $\varepsilon_m$, $f_m$ and $g_m$ in the following way.

CASE 1.

(a) $0 < B(V, \mathcal{F}, \mathcal{G}) < \infty$ and $0 < \varepsilon(V, \mathcal{F}, \mathcal{G}) < \infty$. In this case let $V_m = V$, $l(m) = m$ and define $\varepsilon_m = \varepsilon(V_m, \mathcal{F}_{l(m)}, \mathcal{G}_{l(m)})$. Note that for large $m$, $\varepsilon_m > 0$. Moreover, since both $\mathcal{F}_m$ and $\mathcal{G}_m$ are contained in $B(m)$ it follows that $\varepsilon_m < 2m$. Since $\mathcal{F}_{l(m)}$ and $\mathcal{G}_{l(m)}$ are closed and norm bounded it follows from Lemma 2 of Donoho (1994) that the modulus $\omega(\varepsilon_m, \mathcal{F}_{l(m)}, \mathcal{G}_{l(m)})$ is attained by a pair $f_m \in \mathcal{F}_{l(m)}$ and $g_m \in \mathcal{G}_{l(m)}$.

(b) $\varepsilon(V, \mathcal{F}, \mathcal{G}) = \infty$. In this case let $V_m > V$ be chosen where $V_m \downarrow V$. Then it follows from Lemma 1 that $B(V_m, \mathcal{F}, \mathcal{G}) \to B(V, \mathcal{F}, \mathcal{G})$. So for large $m$, $0 < B(V_m, \mathcal{F}, \mathcal{G}) < \infty$. Now choose an increasing sequence $l(m) \to \infty$ so that $B(V_m, \mathcal{F}_{l(m)}, \mathcal{G}_{l(m)}) > 0$. Now define $\varepsilon_m = \varepsilon(V_m, \mathcal{F}_{l(m)}, \mathcal{G}_{l(m)})$ and once again note that for large $m$, $0 < \varepsilon_m < 2m$. Again $\mathcal{F}_{l(m)}$ and $\mathcal{G}_{l(m)}$ are closed and norm bounded so the modulus $\omega(\varepsilon_m, \mathcal{F}_{l(m)}, \mathcal{G}_{l(m)})$ is attained by a pair $f_m \in \mathcal{F}_{l(m)}$ and $g_m \in \mathcal{G}_{l(m)}$.

CASE 2.

(a) $B(V, \mathcal{F}, \mathcal{G}) = 0$, $B(V', \mathcal{F}, \mathcal{G}) > 0$ for all $V' < V$ and $\omega(\varepsilon, \mathcal{F}, \mathcal{G}) = \sqrt{nV}\varepsilon$ on some interval $0 \leq \varepsilon \leq \varepsilon_0$ for some $\varepsilon_0 > 0$. Let $l(m) = m$ and note that at least for $m$ sufficiently large $0 < B(0, \mathcal{F}_{l(m)}, \mathcal{G}_{l(m)}) < \infty$ and that since $\mathcal{F}_{l(m)} \subseteq \mathcal{F}$ and $\mathcal{G}_{l(m)} \subseteq \mathcal{G}$ it also follows that $B(V, \mathcal{F}_{l(m)}, \mathcal{G}_{l(m)}) = 0$. Lemma 1 shows that for all sufficiently large $m$ there exists a $V_m < V$ such that

$$0 < B(V_m, \mathcal{F}_{l(m)}, \mathcal{G}_{l(m)}) \leq \frac{1}{m}.$$

Now let $\varepsilon_m = \varepsilon(V_m, \mathcal{F}_{l(m)}, \mathcal{G}_{l(m)})$. Then as before it follows that for large $m$, $0 < \varepsilon_m < 2m$. Now since $\mathcal{F}_{l(m)}$ and $\mathcal{G}_{l(m)}$ are closed and norm bounded, the modulus $\omega(\varepsilon_m, \mathcal{F}_{l(m)}, \mathcal{G}_{l(m)})$ is attained by a pair $f_m \in \mathcal{F}_{l(m)}$ and $g_m \in \mathcal{G}_{l(m)}$.

(b) $B(V, \mathcal{F}, \mathcal{G}) = 0$, $B(V', \mathcal{F}, \mathcal{G}) > 0$ for all $0 \leq V' < V$, $\omega(\varepsilon, \mathcal{F}, \mathcal{G}) < \sqrt{nV}\varepsilon$ whenever $\varepsilon > 0$.

Now let $V_m < V$ be chosen where $V_m \uparrow V$. Note that there exists some $V_0 \geq 0$ such that

$$0 < B(V', \mathcal{F}, \mathcal{G}) < \infty$$

for $V_0 \leq V' < V$. Then for large $m$,

$$0 < B(V_m, \mathcal{F}, \mathcal{G}) < \infty.$$

So there is an increasing sequence $l(m) \to \infty$ such that

$$0 < B(V_m, \mathcal{F}_{l(m)}, \mathcal{G}_{l(m)}) < B(V_m, \mathcal{F}, \mathcal{G}) < \infty.$$

We now define $\varepsilon_m = \varepsilon(V_m, \mathcal{F}_{l(m)}, \mathcal{G}_{l(m)})$. It follows once again that $0 < \varepsilon_m < 2m$ for large $m$. Now since $\mathcal{F}_{l(m)}$ and $\mathcal{G}_{l(m)}$ are closed and norm bounded the modulus $\omega(\varepsilon_m, \mathcal{F}_{l(m)}, \mathcal{G}_{l(m)})$ is attained by a pair $f_m \in \mathcal{F}_{l(m)}$ and $g_m \in \mathcal{G}_{l(m)}$.

For $V_m$, $\varepsilon_m$, $f_m$ and $g_m$ as just defined let $u_m = \frac{g_m - f_m}{\varepsilon_m}$ and let

$$a_m = T\left(\frac{f_m + g_m}{2}\right) - \sqrt{nV_m}\left\langle u_m, \frac{f_m + g_m}{2}\right\rangle.$$

For the white noise with drift model let

(21) $$\hat{T}_m = a_m + \sqrt{nV_m} \int u_m(t)\, dY(t)$$

and for the sequence model let

(22) $$\hat{T}_m = a_m + \sqrt{nV_m} \sum_i u_m(i) Y(i).$$

The estimator $\hat{T}_m$ corresponds to the estimator $\hat{T}_V$ defined in Theorem 2 for $V = V_m$, $\mathcal{F} = \mathcal{F}_{l(m)}$ and $\mathcal{G} = \mathcal{G}_{l(m)}$. In the general case we need to take a limit of the estimators $\hat{T}_m$.

Note that $\|u_m\|_2 = 1$ and so there exists a subsequence which converges weakly to some function $u$ where $\|u\|_2 \leq 1$.

Now let $h \in \mathcal{F} \cap \mathcal{G}$. Then since $\|h\|_2 \leq m_0$ for some $m_0 < \infty$ it follows that $h \in \mathcal{F}_{l(m)} \cap \mathcal{G}_{l(m)}$ for all $m \geq m_1$ where $l(m_1) \geq m_0$.

For $m \geq m_1$ note it follows from Theorem 2 that

$$|E_h \hat{T}_m - Th| \leq B(V_m, \mathcal{F}_{l(m)}, \mathcal{G}_{l(m)}) \leq \frac{1}{m}$$

in case 2(a), and in all other cases

$$|E_h \hat{T}_m - Th| \leq B(V_m, \mathcal{F}_{l(m)}, \mathcal{G}_{l(m)}) \leq B(V_m, \mathcal{F}, \mathcal{G}).$$

Note that $B(V_m, \mathcal{F}, \mathcal{G})$ is bounded since it converges to $B(V, \mathcal{F}, \mathcal{G})$. Note also that

$$E\hat{T}_m = a_m + \sqrt{nV_m} \langle u_m, h \rangle$$

and since the norm of $u_m$ is equal to one it follows that $a_m$ is bounded. Hence there is a subsequence of the subsequence used to define $u$ which converges to some finite $a$. Denote this subsubsequence by $m_k^*$.

For the white noise with drift model let

$$\hat{T}_V = a + \sqrt{nV} \int u(t) \, dY(t) \tag{23}$$

and for the sequence model let

$$\hat{T}_V = a + \sqrt{nV} \sum u(i) Y(i). \tag{24}$$

The following theorem shows that this estimator $\hat{T}_V$ which has been formed as a limit of $\hat{T}_m$ trades bias and variance in the general case.

THEOREM 3. *Suppose $\mathcal{F}$ and $\mathcal{G}$ are convex and closed with nonempty intersection. Then the estimator defined by (23) for the white noise with drift model and (24) for the sequence model satisfies*

$$E(\hat{T}_V - E\hat{T}_V)^2 \leq V \tag{25}$$

*and has biases bounded by*

$$\sup_{f \in \mathcal{F}} E_f \hat{T}_V - Tf \leq B(V, \mathcal{F}, \mathcal{G}) \tag{26}$$

*and*

$$\inf_{g \in \mathcal{G}} E_g \hat{T}_V - Tg \geq -B(V, \mathcal{F}, \mathcal{G}). \tag{27}$$

PROOF. Note that (25) follows immediately from the fact that the norm of $u$ is bounded by 1. We shall only give the proof for (26) since the proofs for the other cases are analogous.

First note that the estimator $\hat{T}_m$ as defined in (21) and (22) satisfies the bounds given in Theorem 2. If $f \in \mathcal{F}_{l(m)}$ then

$$E\hat{T}_m - Tf \leq B(V_m, \mathcal{F}_{l(m)}, \mathcal{G}_{l(m)}).$$

Let $m_k^*$ be the subsubsequence along which $a_m$ and $u_m$ converge to $a$ and $u$, respectively. Now for any $f \in \mathcal{F}$, $f \in \mathcal{F}_{l(m_k^*)}$ for large $k$. So

$$\begin{aligned}
E\hat{T}_V - Tf &\leq \limsup_{k \to \infty} E\hat{T}_{m_k^*} - Tf \\
&\leq \limsup_{k \to \infty} B(V_{m_k^*}, \mathcal{F}_{l(m_k^*)}, \mathcal{G}_{l(m_k^*)}) \\
&\leq \limsup_{k \to \infty} B(V_{m_k^*}, \mathcal{F}, \mathcal{G}) \\
&\leq B(V, \mathcal{F}, \mathcal{G}).
\end{aligned}$$

The last step follows from Lemma 1. □

REMARK. Using the Cramér–Rao inequality arguments found in Low (1995) it can be shown that the linear estimator which attains the bounds in the theorem is in fact unique and must actually attain the inequalities. It then follows that the sequence $u_m$ which was used to define the estimator $\hat{T}_V$ actually converges strongly to $u$ and that the sequence $a_m$ actually converges.

**3. Minimax estimator over a finite union of convex sets.** Let $\mathcal{F} = \bigcup_{i=1}^{k} \mathcal{F}_i$ where for $i = 1, \ldots, k$, $\mathcal{F}_i$ are closed convex spaces with nonempty intersections, that is, $\mathcal{F}_i \cap \mathcal{F}_j \neq \varnothing$ for all $i, j$. Our objective is to construct an estimator which is rate optimal for estimating a linear functional $T$ over the parameter space $\mathcal{F}$. Standard two-point testing arguments as, for example, contained in Donoho and Liu (1991a) or Brown and Low (1996) show that the minimax risk for estimating a linear functional $Tf$ over $\mathcal{F}$ is bounded from below by

$$(28) \qquad \inf_{\hat{T}} \sup_{f \in \mathcal{F}} E(\hat{T} - Tf)^2 \geq \frac{1}{8}\omega^2\left(\frac{1}{\sqrt{n}}, \mathcal{F}\right).$$

Let $\hat{T}_i$ be linear estimators which satisfy

$$(29) \qquad \sup_{f \in \mathcal{F}_i} E(\hat{T}_i - Tf)^2 \leq M^2 \omega^2\left(\frac{1}{\sqrt{n}}, \mathcal{F}_i\right)$$

for some $M > 0$. As mentioned in the introduction if $M \geq 1$ such linear estimators are guaranteed to exist and can be constructed by the recipe given in Donoho (1994). In the following discussion $C$ will denote generic constants whereas $M$ will always refer to the bounds given in (29).

For $i \neq j$, let $V_{i,j} = \omega^2(\frac{1}{\sqrt{n}}, \mathcal{F}_i, \mathcal{F}_j)$. Then it follows from the concavity of the modulus that $B(V_{i,j}, \mathcal{F}_i, \mathcal{F}_j)$ as defined by (4) satisfies

$$B(V_{i,j}, \mathcal{F}_i, \mathcal{F}_j) = 2^{-1} \sup_{\varepsilon > 0} (\omega(\varepsilon, \mathcal{F}_i, \mathcal{F}_j) - \sqrt{nV_{i,j}} \varepsilon)$$

$$= 2^{-1} \sup_{\varepsilon \leq 1/\sqrt{n}} \left( \omega(\varepsilon, \mathcal{F}_i, \mathcal{F}_j) - \sqrt{n} \varepsilon \omega\left(\frac{1}{\sqrt{n}}, \mathcal{F}_i, \mathcal{F}_j\right) \right)$$

$$\leq 2^{-1} \omega\left(\frac{1}{\sqrt{n}}, \mathcal{F}_i, \mathcal{F}_j\right).$$

Hence either $B(V_{i,j}, \mathcal{F}_i, \mathcal{F}_j) = 0$ or $0 < B(V_{i,j}, \mathcal{F}_i, \mathcal{F}_j) \leq 2^{-1}\omega(\frac{1}{\sqrt{n}}, \mathcal{F}_i, \mathcal{F}_j)$. In the first case when $B(V_{i,j}, \mathcal{F}_i, \mathcal{F}_j) = 0$ it follows from the definition of $\varepsilon(V_{i,j}, \mathcal{F}_i, \mathcal{F}_j)$ given for case 2(a) that $\varepsilon(V_{i,j}, \mathcal{F}_i, \mathcal{F}_j) = \frac{1}{\sqrt{n}}$. On the other hand if $0 < B(V_{i,j}, \mathcal{F}_i, \mathcal{F}_j) \leq 2^{-1}\omega(\frac{1}{\sqrt{n}}, \mathcal{F}_i, \mathcal{F}_j)$ then $0 < \varepsilon(V_{i,j}, \mathcal{F}_i, \mathcal{F}_j) \leq \frac{1}{\sqrt{n}}$. Hence we know in both cases that $0 < \varepsilon(V_{i,j}, \mathcal{F}_i, \mathcal{F}_j) \leq \frac{1}{\sqrt{n}}$. It follows that when using $V_{i,j} = \omega^2(\frac{1}{\sqrt{n}}, \mathcal{F}_i, \mathcal{F}_j)$ that we are in either case 1(a) or case 2(a) of Section 2.

For $i \neq j$ let $\hat{T}_{i,j}$ be the estimator defined as in Theorem 2 when (9) is attained where $\mathcal{F} = \mathcal{F}_i$, $\mathcal{G} = \mathcal{F}_j$ and $V = V_{i,j} = \omega^2(\frac{1}{\sqrt{n}}, \mathcal{F}_i, \mathcal{F}_j)$. When (9) is not attained the estimator $\hat{T}_{i,j}$ is defined as in Theorem 3. This linear estimator has variance bounded by $\omega^2(\frac{1}{\sqrt{n}}, \mathcal{F}_i, \mathcal{F}_j)$ and bias which satisfies

$$(30) \qquad -2^{-1}\omega\left(\frac{1}{\sqrt{n}}, \mathcal{F}_i, \mathcal{F}_j\right) \leq \inf_{f \in \mathcal{F}_j} (E(\hat{T}_{i,j}) - Tf)$$

and

$$(31) \qquad \sup_{f \in \mathcal{F}_i} (E(\hat{T}_{i,j}) - Tf) \leq 2^{-1}\omega\left(\frac{1}{\sqrt{n}}, \mathcal{F}_i, \mathcal{F}_j\right).$$

Now based on the linear estimators $\hat{T}_{i,j}$ and the linear estimators $\hat{T}_i$, which satisfy (29), define $\hat{z}^u_{i,j}$, $\hat{z}^l_{i,j}$ and $\hat{z}_{i,j}$ by

$$\hat{z}^u_{i,j} = \frac{\hat{T}_{i,j} - \hat{T}_i}{\omega(\frac{1}{\sqrt{n}}, \mathcal{F}_i, \mathcal{F}_j) + M\omega(\frac{1}{\sqrt{n}}, \mathcal{F}_i)},$$

$$\hat{z}^l_{i,j} = \frac{\hat{T}_i - \hat{T}_{j,i}}{\omega(\frac{1}{\sqrt{n}}, \mathcal{F}_j, \mathcal{F}_i) + M\omega(\frac{1}{\sqrt{n}}, \mathcal{F}_i)}$$

and

$$\hat{z}_{i,j} = \max(\hat{z}^u_{i,j}, \hat{z}^l_{i,j}).$$

Note that $\hat{z}^u_{i,j}$ and $\hat{z}^l_{i,j}$ are normally distributed and satisfy

(32) $$\max(\text{Var}(\hat{z}^u_{i,j}), \text{Var}(\hat{z}^l_{i,j})) \leq 1.$$

Finally define the estimator of the linear functional $T$ as

(33) $$\hat{T}^* = \hat{T}_{\hat{i}} \quad \text{with } \hat{i} = \arg\min_i \left(\sup_{j \neq i} \hat{z}_{i,j}\right).$$

The analysis of the mean squared error of $\hat{T}^*$ is facilitated by the following lemma which bounds the probability that $\hat{T}^* = \hat{T}_j$ when the magnitude of the bias of $\hat{T}_j$ is large.

LEMMA 2. *Suppose $f \in \mathcal{F}_i$ and for some $j \neq i$, $|E\hat{T}_j - Tf| \geq \gamma M\omega(\frac{1}{\sqrt{n}}, \mathcal{F})$ where $\gamma \geq 3$. Then*

$$P(\hat{i} = j) \leq 2k \exp\left(-\frac{(\gamma - 3)^2}{32}\right).$$

PROOF. Note that if $f \in \mathcal{F}_i$ and $j \neq i$, then from (29) and (31),

(34) $$E\hat{z}^u_{i,j} = \frac{E(\hat{T}_{i,j} - Tf) - E(\hat{T}_i - Tf)}{\omega(\frac{1}{\sqrt{n}}, \mathcal{F}_i, \mathcal{F}_j) + M\omega(\frac{1}{\sqrt{n}}, \mathcal{F}_i)} \leq 1$$

and from (29) and (30),

(35) $$E\hat{z}^l_{i,j} = \frac{E(\hat{T}_i - Tf) - E(\hat{T}_{j,i} - Tf)}{\omega(\frac{1}{\sqrt{n}}, \mathcal{F}_j, \mathcal{F}_i) + M\omega(\frac{1}{\sqrt{n}}, \mathcal{F}_i)} \leq 1.$$

Now suppose that $f \in \mathcal{F}_i$. We shall only give details of the proof when

$$E\hat{T}_j - Tf \geq \gamma M\omega\left(\frac{1}{\sqrt{n}}, \mathcal{F}\right)$$

as the case when

$$E\hat{T}_j - Tf \leq -\gamma M\omega\left(\frac{1}{\sqrt{n}}, \mathcal{F}\right)$$

is handled in a similar way. When $E\hat{T}_j - Tf \geq \gamma M\omega(\frac{1}{\sqrt{n}}, \mathcal{F})$ then it follows from (31) that

(36) $$\begin{aligned} E\hat{z}^l_{j,i} &= \frac{E(\hat{T}_j - Tf) - E(\hat{T}_{i,j} - Tf)}{\omega(\frac{1}{\sqrt{n}}, \mathcal{F}_i, \mathcal{F}_j) + M\omega(\frac{1}{\sqrt{n}}, \mathcal{F}_j)} \\ &\geq \frac{\gamma M\omega(\frac{1}{\sqrt{n}}, \mathcal{F}) - \frac{1}{2}\omega(\frac{1}{\sqrt{n}}, \mathcal{F}_i, \mathcal{F}_j)}{\omega(\frac{1}{\sqrt{n}}, \mathcal{F}_i, \mathcal{F}_j) + M\omega(\frac{1}{\sqrt{n}}, \mathcal{F}_j)} \\ &\geq \frac{\gamma - 1}{2}. \end{aligned}$$

Now without loss of generality suppose that $i = 1$ and that $j = 2$. Then if $\hat{i} = 2$ note that $\hat{z}_{2,1} \leq \sup_{r \neq 1}(\hat{z}_{1,r})$ and since $\hat{z}_{2,1}^l \leq \hat{z}_{2,1}$ it follows that

$$P(\hat{i} = 2) \leq \sum_{r=2}^{k} P(\hat{z}_{2,1}^l - \hat{z}_{1,r} \leq 0)$$
$$\leq \sum_{r=2}^{k} \{P(\hat{z}_{2,1}^l - \hat{z}_{1,r}^l \leq 0) + P(\hat{z}_{2,1}^l - \hat{z}_{1,r}^u \leq 0)\}.$$

Now by (32) $\hat{z}_{2,1}^l - \hat{z}_{1,r}^l$ and $\hat{z}_{2,1}^l - \hat{z}_{1,r}^u$ both have normal distributions with variance less than or equal to 4 and by (34)–(36) means greater than or equal to $\frac{\gamma - 3}{2}$ and the lemma now follows from the bound on a standard normal random variable $Z$,

$$P(Z > t) \leq \exp\left(-\frac{t^2}{2}\right)$$

which holds for all $t \geq 0$.  □

Although our main focus is on mean squared error we shall consider the more general case of $p$th power loss. Such general cases are important in the theory of adaptation [see Cai and Low (2002)]. Lemma 2 can be used to bound the risk of the estimator $\hat{T}^*$ defined by (33) as in the following theorem.

THEOREM 4. *Suppose either the white noise model (1) or the sequence model (2) is given. Let $\mathcal{F} = \bigcup_{i=1}^{k} \mathcal{F}_i$ where $k \geq 2$ and $\mathcal{F}_i$ are closed convex sets with $\mathcal{F}_i \cap \mathcal{F}_j \neq \varnothing$ for all $i, j$. Let $\hat{T}^*$ be the estimator of the linear functional $T$ defined as in (33). Then for $p \geq 1$,*

$$(37) \quad \sup_{f \in \mathcal{F}} E|\hat{T}^* - Tf|^p \leq C(p) M^p (\ln k)^{p/2} \omega^p\left(\frac{1}{\sqrt{n}}, \mathcal{F}\right)$$

*where the constant $C(p)$ is independent of $M$, $k$ and $n$.*

REMARK. Note that we can always find linear estimators $\hat{T}_i$ for which $M \leq 1$ in (29) and so the theorem yields an upper bound on the minimax risk over $\mathcal{F}$ which only depends on the modulus and the number $k$. There is also a minimax lower bound for the $p$th power loss analogous to that given in (28),

$$(38) \quad \inf_{\hat{T}} \sup_{f \in \mathcal{F}} E|\hat{T} - Tf|^p \geq b(p) \omega^p\left(\frac{1}{\sqrt{n}}, \mathcal{F}\right)$$

for some constant $b(p) > 0$. By comparing the upper bound in (37) with this bound it is clear that for fixed finite $k$ the estimator $\hat{T}^*$ is rate optimal over $\mathcal{F}$. It is also worth noting that sometimes the lower bound is not asymptotically sharp when $k$ is finite but grows with $n$. In Section 5 we give examples where $k$ grows with $n$ and the optimal rate is given by the upper bound in equation (37).

In the theory of adaptation the goal is to find a procedure which is simultaneously near minimax over a collection of parameter spaces. If a collection of convex parameter spaces is not nested then the largest of the minimax risks for each convex parameter space may be smaller than the minimax risk over the union of the convex parameter spaces [see, e.g., Cai and Low (2002)]. In such cases an appropriate benchmark for the maximum risk of an adaptive estimator is given by the bound in Theorem 4.

The proof of Theorem 4 is based on Lemma 2 and the following bound on the tail probabilities of a maximum of Gaussian random variables.

LEMMA 3. *Let $X_i$, $i = 1, \ldots, m$, be normal random variables with means $\mu_i$ and standard deviations $\sigma_i \leq \sigma$. Suppose that $|\mu_i - \mu| \leq \gamma$ for $i = 1, \ldots, m$, and $c > 0$ is a constant. Then*

$$(39) \qquad P\left(\max_{1 \leq i \leq m} |X_i - \mu| \geq \gamma + \sqrt{c \ln m}\, \sigma\right) \leq m^{1-c/2}.$$

PROOF. We shall assume that $m \geq 2$ and that $c \geq 2$, since otherwise the bound is trivial. Denote by $Z$ a standard Gaussian random variable. Then

$$P\left(\max_{1 \leq i \leq m} |X_i - \mu| \geq \gamma + \sqrt{c \ln m}\, \sigma\right) \leq \sum_{i=1}^{m} P(|X_i - \mu| \geq \gamma + \sqrt{c \ln m}\, \sigma)$$

$$\leq m P(|Z| \geq \sqrt{c \ln m})$$

$$\leq m^{1-c/2}.$$

The last inequality follows from standard bounds on tail probabilities of Gaussian distributions once we note that $c \ln m \geq 1$ when $c \geq 2$ and $m \geq 2$. □

PROOF OF THEOREM 4. Let $\lambda \geq 1$ and $D_\lambda = 3 + \sqrt{32 \lambda \ln 2k^3}$ and we will write $D_p$ when $\lambda = p$. Then it is easy to check from Lemma 2 that if $|E\hat{T}_i - Tf| \geq D_\lambda M \omega(\frac{1}{\sqrt{n}}, \mathcal{F})$ then $P(\hat{i} = i) \leq \frac{1}{k^{3\lambda - 1}}$.

Let

$$I_1 = \left\{i : |E\hat{T}_i - Tf| \geq D_p M \omega\left(\frac{1}{\sqrt{n}}, \mathcal{F}\right)\right\},$$

$$I_2 = \left\{i : |E\hat{T}_i - Tf| < D_p M \omega\left(\frac{1}{\sqrt{n}}, \mathcal{F}\right)\right\}.$$

We then have

$$
\begin{aligned}
E|\hat{T}^* - Tf|^p &= \sum_{i \in I_1} E(|\hat{T}_i - Tf|^p \mathbb{1}(\hat{i} = i)) + \sum_{i \in I_2} E(|\hat{T}_i - Tf|^p \mathbb{1}(\hat{i} = i)) \\
&\leq \sum_{i \in I_1} (E(\hat{T}_i - Tf)^{2p})^{1/2} (P(\hat{i} = i))^{1/2} + E\left(\max_{i \in I_2} |\hat{T}_i - Tf|^p\right).
\end{aligned}
\tag{40}
$$

Now note that, if $X$ has a normal distribution with mean $\mu$ and variance $\sigma^2$ and $p \geq 1$, then a straightforward calculation shows that

$$(EX^{2p})^{1/2} \leq 2^p(|\mu^p| + a_{2p}^{1/2}\sigma^p)$$

where $a_j = E|Z|^j$ for $Z$ a standard Gaussian random variable. If $i \in I_1$, then for some $\lambda \geq p$,

$$|E\hat{T}_i - Tf| = D_\lambda M \omega\left(\frac{1}{\sqrt{n}}, \mathcal{F}\right),$$

and so for $i \in I_1$ and such a choice of $\lambda \geq p$,

$$(E|\hat{T}_i - Tf|^{2p})^{1/2}(P(\hat{i} = i))^{1/2} \leq 2^p M^p (a_{2p}^{1/2} + D_\lambda^p) \omega^p\left(\frac{1}{\sqrt{n}}, \mathcal{F}\right) \frac{1}{k^{(3\lambda-1)/2}}.$$

Now note that, if $k \geq 2$,

$$\sup_{\lambda \geq p} \frac{D_\lambda^p}{k^{(3\lambda-1)/2}} = \frac{D_p^p}{k^{(3p-1)/2}} \leq \frac{D_p^p}{k},$$

and hence

$$\sum_{i \in I_1} (E|\hat{T}_i - Tf|^{2p})^{1/2}(P(\hat{i} = i))^{1/2} \leq 2^p M^p (a_{2p}^{1/2} + D_p^p) \omega^p\left(\frac{1}{\sqrt{n}}, \mathcal{F}\right). \tag{41}$$

Let $m$ be the cardinality of the set $I_2$. Now note that, if $m \leq 1$, then

$$E\left(\max_{i \in I_2} |\hat{T}_i - Tf|^p\right) \leq B(p) \omega^p\left(\frac{1}{\sqrt{n}}, \mathcal{F}\right) \tag{42}$$

and the theorem now follows from (40)–(42). If $m \geq 2$, then

$$
\begin{aligned}
E&\left(\max_{i \in I_2} |\hat{T}_i - Tf|^p\right) \\
&\leq (D_p + \sqrt{3 \ln m})^p M^p \omega^p\left(\frac{1}{\sqrt{n}}, \mathcal{F}\right) \\
&\quad + M^p \omega^p\left(\frac{1}{\sqrt{n}}, \mathcal{F}\right) \\
&\quad \times \sum_{l=3}^{\infty} \Big\{(D_p + \sqrt{l \ln m})^p
\end{aligned}
$$

$$\times P\Big((D_p + \sqrt{(l-1)\ln m}\,)M\omega\Big(\frac{1}{\sqrt{n}},\mathcal{F}\Big) \leq \max_{i\in I_2}|\hat{T}_i - Tf|$$
$$\leq (D_p + \sqrt{l\ln m}\,)M\omega\Big(\frac{1}{\sqrt{n}},\mathcal{F}\Big)\Big)\Big\}$$
$$\leq (D_p + \sqrt{3\ln m}\,)^p M^p \omega^p\Big(\frac{1}{\sqrt{n}},\mathcal{F}\Big)$$
$$+ M^p \omega^p\Big(\frac{1}{\sqrt{n}},\mathcal{F}\Big)$$
$$\times \sum_{l=3}^{\infty}\Big\{(D_p + \sqrt{l\ln m}\,)^p$$
$$\times P\Big(\max_{i\in I_2}|\hat{T}_i - Tf| \geq (D_p + \sqrt{(l-1)\ln m}\,)M\omega\Big(\frac{1}{\sqrt{n}},\mathcal{F}\Big)\Big)\Big\}.$$

Note that it follows from the definition of $I_2$ and the fact that the variance of $\hat{T}_i$ is bounded by $M^2\omega^2(\frac{1}{\sqrt{n}},\mathcal{F}))$ and Lemma 3 that

$$(43)\quad P\Big(\max_{i\in I_2}|\hat{T}_i - Tf| \geq (D_p + \sqrt{(l-1)\ln m}\,)M\omega\Big(\frac{1}{\sqrt{n}},\mathcal{F}\Big)\Big) \leq m^{-(l-3)/2}.$$

Hence

$$E\Big(\max_{i\in I_2}|\hat{T}_i - Tf|^p\Big)$$
$$(44)\quad \leq \Big[(D_p + \sqrt{3\ln m}\,)^p + \sum_{l=3}^{\infty}(D_p + \sqrt{l\ln m}\,)^p m^{-(l-3)/2}\Big]M^p\omega^p\Big(\frac{1}{\sqrt{n}},\mathcal{F}\Big)$$
$$\leq B(p)(\ln k)^{p/2} M^p \omega^p\Big(\frac{1}{\sqrt{n}},\mathcal{F}\Big).$$

The theorem now follows on combining (40), (41) and (44).  □

**4. Linear estimators.** We now consider the performance of linear procedures. As mentioned in the Introduction, the optimal linear procedure is within a small constant factor of the minimax risk when the parameter space is convex. The following theorem considers the case when the parameter space is nonconvex. Let $\mathcal{F}$ denote a parameter set and let C.Hull($\mathcal{F}$) denote the convex hull of $\mathcal{F}$.

THEOREM 5. *Consider the white noise model* (1) *or the sequence model* (2). *The minimax linear risk over a parameter set $\mathcal{F}$ is the same as the minimax linear risk over the convex hull of $\mathcal{F}$, that is,*

$$R_A^*(n;\mathcal{F}) = R_A^*(n;\text{C.Hull}(\mathcal{F})).$$

This theorem is a direct consequence of the following result.

THEOREM 6. *Let $\hat{T}$ be a linear estimator of $Tf$ where $T$ is a linear functional. Then for any $\mathcal{F}$*

$$\sup_{f \in \mathcal{F}} E_f(\hat{T} - Tf)^2 = \sup_{f \in \mathrm{C.Hull}(\mathcal{F})} E_f(\hat{T} - Tf)^2. \tag{45}$$

PROOF. Since $\mathcal{F} \subseteq \mathrm{C.Hull}(\mathcal{F})$, it is obvious that

$$\sup_{f \in \mathcal{F}} E(\hat{T} - Tf)^2 \leq \sup_{f \in \mathrm{C.Hull}(\mathcal{F})} E(\hat{T} - Tf)^2.$$

Let $f \in \mathrm{C.Hull}(\mathcal{F})$ and $f = \sum_i \lambda_i f_i$ with $f_i \in \mathcal{F}$, $\lambda_i \geq 0$ and $\sum_i \lambda_i = 1$. Then the squared bias

$$(E_f \hat{T} - Tf)^2 = \left( \sum_i \lambda_i (E_{f_i} \hat{T} - Tf_i) \right)^2 \leq \left( \sum_i \lambda_i |E_{f_i} \hat{T} - Tf_i| \right)^2$$

$$\leq \max_i |E_{f_i} \hat{T} - Tf_i|^2 \leq \sup_{f \in \mathcal{F}} (E_f \hat{T} - Tf)^2.$$

It then follows from the fact that a linear estimator has constant variance that

$$\sup_{f \in \mathrm{C.Hull}(\mathcal{F})} E(\hat{T} - Tf)^2 \leq \sup_{f \in \mathcal{F}} E(\hat{T} - Tf)^2. \qquad \square$$

Note that equation (45) is not necessarily true for nonlinear procedures. The following corollary is a direct consequence of Theorem 5.

COROLLARY 1.

$$R_A^*(n; \mathcal{F}) = \sup_{\varepsilon > 0} \omega^2(n, \mathrm{C.Hull}(\mathcal{F})) \frac{1/(4n)}{1/n + \varepsilon^2/4} \tag{46}$$

*and*

$$\frac{1}{8} \omega^2 \left( \frac{2}{\sqrt{n}}, \mathrm{C.Hull}(\mathcal{F}) \right) \leq R_A^*(n; \mathcal{F}) \leq \frac{1}{4} \omega^2 \left( \frac{2}{\sqrt{n}}, \mathrm{C.Hull}(\mathcal{F}) \right). \tag{47}$$

Thus the minimax linear risk is determined by the modulus of continuity over the convex hull of $\mathcal{F}$, not over $\mathcal{F}$ itself. In the case that $\omega(\varepsilon, \mathrm{C.Hull}(\mathcal{F})) \gg \omega(\varepsilon, \mathcal{F})$, linear procedures will perform poorly. Examples which illustrate this point are contained in the next section.

**5. Examples.** In this section we discuss examples where the modulus of continuity over the convex hull of the parameter space is much larger than the modulus of continuity over the parameter space. Since the performance of the optimal linear procedure is determined by the modulus of the convex hull of the parameter space linear procedures perform badly in these cases. On the other hand, the nonlinear procedure introduced in Section 3 is within a constant factor of the minimax risk.

5.1. *Estimating functions at a point.* Suppose we observe the white noise model (1) over the interval $[-\frac{1}{2}, \frac{1}{2}]$ and we wish to estimate $Tf = f(0)$.

We recall that a function is Lip($\alpha$) ($0 < \alpha \leq 1$) over an interval $[a, b]$ if

$$|f(x) - f(y)| \leq |x-y|^\alpha \qquad \text{for all } x, y \in [a,b].$$

Let

$$\mathcal{F}_1 = \{f : f \text{ is continuous on } [-\tfrac{1}{2}, \tfrac{1}{2}]$$
$$\text{with maximum at } 0 \text{ and } f \text{ is Lip}(1) \text{ over } [-\tfrac{1}{2}, 0]\}$$

and

$$\mathcal{F}_2 = \{f : f \text{ is continuous on } [-\tfrac{1}{2}, \tfrac{1}{2}]$$
$$\text{with maximum at } 0 \text{ and } f \text{ is Lip}(\tfrac{1}{2}) \text{ over } [0, \tfrac{1}{2}]\}.$$

Let $\mathcal{F} = \mathcal{F}_1 \cup \mathcal{F}_2$. The parameter spaces $\mathcal{F}_1$ and $\mathcal{F}_2$ are both convex, but $\mathcal{F}$ is nonconvex. It is easy to see that

$$\text{C.Hull}(\mathcal{F}) = \{\text{All continuous functions over } [-\tfrac{1}{2}, \tfrac{1}{2}] \text{ with maximum at } 0\}.$$

The convex hull of $\mathcal{F}$ is "much larger" than $\mathcal{F}$. By straightforward calculations it is easy to verify that for $Tf = f(0)$ and small $\varepsilon > 0$,

$$\omega(\varepsilon, \mathcal{F}_1) = \omega(\varepsilon, \mathcal{F}_2, \mathcal{F}_1) = 3^{1/3}\varepsilon^{2/3},$$
$$\omega(\varepsilon, \mathcal{F}_2) = \omega(\varepsilon, \mathcal{F}_1, \mathcal{F}_2) = 2^{1/4}\varepsilon^{1/2}(1 + o(1))$$

so $\omega(\varepsilon, \mathcal{F}) = 2^{1/4}\varepsilon^{1/2}(1 + o(1))$. But $\omega(\varepsilon, \text{C.Hull}(\mathcal{F})) = \infty$.

It follows from Theorem 4 that the minimax mean squared error rate of convergence for estimating the linear functional $Tf = f(0)$ is $n^{-1/2}$. However, the maximum risk of any linear estimator over $\mathcal{F}$ is not even bounded. [This follows from the fact that $\omega(\varepsilon, \text{C.Hull}(\mathcal{F})) = \infty$.] In other words, linear estimators do not work at all in this case.

5.2. *Estimating a linear functional of nearly black objects.* In this example we consider the Gaussian sequence model

(48) $$y_i = f_i + n^{-1/2} z_i, \qquad i = 1, \ldots, n,$$

where $z_i \stackrel{\text{i.i.d.}}{\sim} N(0,1)$. The size of the vector, $n$, is assumed large; we are interested in asymptotics in which the number of variables is large. We assume that the vector $f$ is sparse: only a small fraction of components are nonzero, and the indices, or locations of the nonzero components are not known in advance.

Denote the $\ell_0$ quasi-norm by $\|f\|_0 = \text{Card}(\{i : f_i \neq 0\})$. Fix $k_n$, the collection of vectors with at most $k_n$ nonzero entries is

$$\mathcal{F} = \ell_0(k_n) = \{f \in \mathbb{R}^n : \|f\|_0 \leq k_n\}.$$

Following Donoho, Johnstone, Hoch and Stern (1992), we call a setting nearly black when the fraction of nonzero components $k_n/n \approx 0$, by analogy with night-sky images. In this example we assume that $k_n$ is known and $k_n \leq C n^\varepsilon$ where $\varepsilon < 1/2$.

A motivation for this model is provided by wavelet analysis, since the wavelet representation of many smooth and piecewise smooth signals is sparse and nearly black in this sense [see, e.g., Donoho, Johnstone, Kerkyacharian and Picard (1995)]. For estimating the whole object, this model has also been studied in Donoho, Johnstone, Hoch and Stern (1992) and Abramovich, Benjamini, Donoho and Johnstone (2000).

In the present paper we are interested in estimating the linear functional of the unknown vector $f$ given by

$$Tf = \sum_{i=1}^n f_i.$$

Let $\mathcal{I}(k,n)$ be the class of all subsets of $\{1,\ldots,n\}$ of $k$ elements and for $I \in \mathcal{I}(k,n)$ let

$$\mathcal{F}_I = \{f \in \mathbb{R}^n : f_j = 0 \ \forall j \notin I\}.$$

Note that $\mathcal{F}_I$ is a $k_n$-dimensional subspace spanned by the coordinates in $I$. These are obviously convex and $\mathcal{F} = \cup \mathcal{F}_I$ where the union is taken over $I$ in the set $\mathcal{I}(k_n,n)$. From now on we shall assume that $I$ is in the set $\mathcal{I}(k_n,n)$.

Linear procedures perform poorly over $\mathcal{F}$. In fact it is easy to see that the convex hull of $\mathcal{F}$ is the whole of $\mathbb{R}^n$ and

$$\omega(\varepsilon, \text{C.Hull}(\mathcal{F})) = \omega(\varepsilon, \mathbb{R}^n) = \sqrt{n}\varepsilon.$$

It then follows from Theorem 5 that any linear estimator must have maximum mean squared error over $\mathcal{F}$ of at least 1. In fact it is easy to see that the best linear procedure is simply $\hat{T} = \sum_{i=1}^n y_i$.

Nonlinear procedures can perform much better. Our general construction given in Section 3 starts with linear estimators constructed assuming that $f \in \mathcal{F}_I$. In this case it is natural to start with $\hat{T}_I$ the minimax estimator over $\mathcal{F}_I$ since this estimator is linear, unbiased over $\mathcal{F}_I$ and has variance equal to $\frac{k_n}{n}$. $\hat{T}_I$ is in fact just the sum of $y_j$ with $j \in I$. In this example equation (29) holds for all $I \in \mathcal{I}(k_n, n)$ with $M = 1$ and $\omega^2(\frac{1}{\sqrt{n}}, \mathcal{F}_I) = \frac{k_n}{n}$.

The construction of $\hat{T}^*$ is also based on the modulus $\omega(\varepsilon, \mathcal{F}_I, \mathcal{F}_J)$ between $\mathcal{F}_I$ and $\mathcal{F}_J$. Note that a least favorable pair of parameters is given by one parameter which has the $k_n$ coefficients in $J$ all equal to some given value $a > 0$ and the rest zero and the second parameter has the coefficients in $J \setminus I$ equal to $-a$ and the rest zero. By choosing $a$ so that the $l_2$ distance between these parameters is equal to $\varepsilon$ it is easy to check that

$$\omega(\varepsilon, \mathcal{F}_I, \mathcal{F}_J) = \sqrt{\mathrm{Card}(I \cup J)}\varepsilon$$

and consequently

$$\omega(\varepsilon, \mathcal{F}) = \sqrt{2k_n}\varepsilon.$$

Now let $\hat{T}_{I,J}$ be defined as in Section 3. It is easy to see in this case that

$$\hat{T}_{I,J} = \sum_{l \in I \cup J} y_l.$$

Let $N$ be the number of parameter spaces. Then $N$ is equal to $n$ choose $k_n$ and it is easy to see that

$$N = \binom{n}{k_n} \leq n^{k_n}.$$

It then follows from Theorem 4 that if $\hat{T}^*$ is defined by (33), then

(49) $$\sup_{f \in \mathcal{F}} E(\hat{T}^* - Tf)^2 \leq C \frac{k_n^2 \ln n}{n}.$$

The following theorem shows that the estimator $\hat{T}^*$ is in fact rate optimal. The theorem gives a minimax lower bound based on using a mixture prior and a constrained risk inequality introduced in Brown and Low (1996).

THEOREM 7. *Let $Tf = \sum_{i=1}^{n} f_i$. Suppose that $n \geq 4$ and that $k_n < n^\varepsilon$ with $\varepsilon < 1/2$. Then*

(50) $$\inf_{\hat{T}} \sup_{f \in \mathcal{F}} E(\hat{T} - Tf)^2 \geq \frac{1}{121} \frac{k_n^2}{n} \ln\left(\frac{n}{k_n^2}\right).$$

REMARK. Comparing the minimax lower bound (50) with the risk upper bound for $\hat{T}^*$, for $k_n < Cn^\varepsilon$ with $\varepsilon < 1/2$, the estimator $\hat{T}^*$ is within a constant factor of the minimax risk. For example, for $k_n = n^\varepsilon$ with $\varepsilon < 1/2$, the risk of $\hat{T}^*$ converges at the rate of $n^{-(1-2\varepsilon)} \log n$ which is optimal.

PROOF OF THEOREM 7. In the proof we will omit the subscript in $k_n$ and simply write $k$ for $k_n$. Let $\psi_\mu$ be the density of a normal distribution with mean $\mu$ and variance $\frac{1}{n}$. And for $I \in \mathcal{I}(k,n)$ let

$$g_I(y_1, \ldots, y_n) = \prod_{j=1}^{n} \psi_{f_j}(y_j)$$

where $f_j = \frac{\rho}{\sqrt{n}} \mathbb{1}(j \in I)$. Finally let

$$g = \frac{1}{\binom{n}{k}} \sum_{I \in \mathcal{I}(k,n)} g_I$$

and $f = \prod_{j=1}^{n} f_0$ be the density of $n$ independent normal random variables each with mean 0 and variance $\frac{1}{n}$. Note that a similar mixture prior was used in Baraud (2000) to give lower bounds in a nonparametric testing problem. Now note that if

$$E_{g_I}\left(\delta - k\frac{\rho}{\sqrt{n}}\right)^2 \leq C$$

for all $I \in \mathcal{I}(k,n)$ then it follows that

$$E_g\left(\delta - k\frac{\rho}{\sqrt{n}}\right)^2 \leq C.$$

We will now apply the constrained risk inequality of Brown and Low (1996). First we need to calculate a chi-squared distance between $f$ and $g$. This is done as follows. Note that

$$\int \frac{g^2}{f} = \frac{1}{\binom{n}{k}^2} \sum_{I \in \mathcal{I}(k,n)} \sum_{I' \in \mathcal{I}(k,n)} \int \frac{g_I g_{I'}}{f}$$

and simple calculations show that

$$\int \frac{g_I g_{I'}}{f} = \exp(j\rho^2)$$

where $j$ is the number of points in the set $I \cap I'$. It follows that

$$\int \frac{g^2}{f} = E \exp(J\rho^2)$$

where $J$ has a hypergeometric distribution

$$P(J = j) = \frac{\binom{k}{j}\binom{n-k}{k-j}}{\binom{n}{k}}.$$

Now note that from Feller [(1968), page 59],

$$P(J=j) \leq \binom{k}{j}\left(\frac{k}{n}\right)^j\left(1-\frac{k}{n}\right)^{k-j}\left(1-\frac{k}{n}\right)^{-k}.$$

Now suppose that $n \geq 4$ and that $k < n^{1/2}$. Then

$$\left(1-\frac{k}{n}\right)^{-k} \leq 4^{k^2/n} \leq 4$$

and hence

$$P(J=j) \leq 4\binom{k}{j}\left(\frac{k}{n}\right)^j\left(1-\frac{k}{n}\right)^{k-j}.$$

It now follows that if $n \geq 4$ and $k < n^{1/2}$ then

$$\int \frac{g^2}{f} = E\exp(J\rho^2)$$

$$\leq 4\left(1-\frac{k}{n}+\frac{k}{n}e^{\rho^2}\right)^k.$$

Now take $\rho = \sqrt{\ln \frac{n}{k^2}}$ and it follows that

$$\int \frac{g^2}{f} \leq 4\left(1+\frac{1}{k}\right)^k$$

$$\leq 4e.$$

It then follows from the constrained risk inequality in Brown and Low (1996) that if

$$(51) \qquad E_f(\delta-0)^2 \leq c_1 \frac{k^2}{n}\ln\frac{n}{k^2}$$

then

$$(52) \qquad \begin{aligned}E_g\left(\delta-k\frac{\rho}{\sqrt{n}}\right)^2 &\geq \frac{k^2}{n}\ln\frac{n}{k^2} - 4e\frac{k}{\sqrt{n}}\sqrt{\ln\frac{n}{k^2}}\sqrt{c_1\frac{k^2}{n}\ln\frac{n}{k^2}}\\ &= \frac{k^2}{n}\ln\frac{n}{k^2}(1-4e\sqrt{c_1}).\end{aligned}$$

The theorem now follows on taking $c_1 = 1 + 8e^2 - 4e\sqrt{1+4e^2}$. $\square$

5.3. *Structured nearly black objects.* We will now consider an example under the Gaussian sequence model (48) where most of the coordinates are zero but where we shall also assume that the $k_n$ nonzero coordinates appear consecutively and that $0 \leq k_n \leq n^\varepsilon$ for some $\varepsilon < 1$. Again $k_n$ is assumed to be known. Let

$$\mathcal{F}(a, k_n) = \{f \in R^n : f_i = 0 \text{ unless } a \leq i \leq a+k_n-1\}$$

and

$$\mathcal{F} = \bigcup_{a=1}^{n-k_n} \mathcal{F}(a, k_n).$$

We call members of $\mathcal{F}$ structured nearly black objects. It is easy to see that the convex hull of $\mathcal{F}$ is again the whole of $\mathbb{R}^n$. It thus follows from Theorem 5 that linear procedures perform poorly for estimating $Tf$ over $\mathcal{F}$.

Let $\hat{T}_a = \sum_{i=a}^{a+k_n-1} y_i$. Then $\hat{T}_{a,b}$ as defined in Section 3 is given by

$$\hat{T}_{a,b} = \sum y_i \mathbb{1}(i \in [a, a+k_n-1] \cup [b, b+k_n-1]).$$

Note that $\mathcal{F}$ is a union of only $n - k_n$ convex sets and so it then follows from Theorem 4 that if $\hat{T}^*$ is now defined by (33) then

(53) $$\sup_{f \in \mathcal{F}} E(\hat{T}^* - Tf)^2 \leq C \frac{k_n \ln n}{n}.$$

Equation (53) gives an upper bound for the minimax risk. We shall now show that this upper bound is rate sharp. In fact we shall show that if $n \geq 4$ and $k_n < n^\varepsilon$ with $\varepsilon < 1$, then

(54) $$\inf_{\hat{T}} \sup_{f \in \mathcal{F}} E(\hat{T} - Tf)^2 \geq \frac{1}{18} \frac{k_n}{n} \ln\left(\frac{3n}{k_n}\right).$$

This can be seen as follows. Denote the index sets $I_a = \{i : a \leq i \leq a+k_n-1\}$ and let $\mathcal{I}(k_n, n) = \bigcup_{a=1}^{n-k_n} I_a$. As in the previous example let $\psi_f$ be the density of a normal distribution with mean $f$ and variance $\frac{1}{n}$. And for $I \in \mathcal{I}(k_n, n)$ let

$$g_I(y_1, \ldots, y_n) = \prod_{j=1}^{n} \psi_{f_j}(y_j)$$

where $f_j = \frac{\rho}{\sqrt{n}} \mathbb{1}(j \in I)$. Finally let $g = \frac{1}{n-k_n} \sum_{I=1}^{n-k_n} g_I$ and $f = \prod_{j=1}^{n} f_0$ be the density of $n$ independent normal random variables each with mean 0 and variance $\frac{1}{n}$. Following the argument in the previous example we note that

$$\int \frac{g^2}{f} = E \exp(J\rho^2)$$

where this time $J$ satisfies

$$P(J = 0) = \frac{n - k_n}{n}$$

and for $1 \leq i \leq k_n$,

$$P(J = i) = \frac{1}{n}.$$

Hence

$$\int \frac{g^2}{f} \leq 1 + \frac{k_n}{n} \exp(k_n \rho^2).$$

Now set

$$\rho = \sqrt{\frac{1}{k_n}} \sqrt{\ln\left(\frac{3n}{k_n}\right)}.$$

Then $\int \frac{g^2}{f} \leq 4$ and (54) now follows as in (51) and (52).

**Acknowledgment.** We thank two referees for their thorough and useful comments which have helped to improve the presentation of the paper.

Department of Statistics  
The Wharton School  
University of Pennsylvania  
Philadelphia, Pennsylvania 19104-6340  
USA  
e-mail: tcai@wharton.upenn.edu  
e-mail: lowm@wharton.upenn.edu